\documentclass{article}
\usepackage{amssymb,amsmath,amsthm,graphicx}

\def\qed{\hfill {\hbox{${\vcenter{\vbox{               
   \hrule height 0.4pt\hbox{\vrule width 0.4pt height 6pt
   \kern5pt\vrule width 0.4pt}\hrule height 0.4pt}}}$}}}

\def\utr{\, \underline{\ast}\, }
\def\otr{\, \overline{\ast}\, }
\def\vtr{\! \circledast\! }

\newtheorem{theorem}{Theorem}

\newtheorem{corollary}[theorem]{Corollary}

\theoremstyle{definition}
\newtheorem{example}{Example}
\newtheorem{definition}{Definition}
\newtheorem{remark}{Remark}

\date{}

\title{\Large \textbf{Twisted Virtual Bikeigebras and Twisted Virtual Handlebody-Knots}}

\author{Sam Nelson\footnote{Email: Sam.Nelson@cmc.edu. Partially supported by Simons Foundation collaboration grant 316709}\and
Yuqi Zhao\footnote{Email: YZhao17@students.claremontmckenna.edu.}} 

\begin{document}
\maketitle

\begin{abstract}
We generalize unoriented handlebody-links to the twisted virtual case, 
obtaining Reidemeister moves for handlebody-links in ambient spaces of the
form $\Sigma\times [0,1]$ for $\Sigma$ a compact closed 2-manifold up to
stable equivalence. We introduce a related algebraic structure known as 
twisted virtual bikeigebras whose axioms are motivated by the twisted virtual
handlebody-link Reidemeister moves. We use twisted virtual bikeigebras to 
define $X$-colorability for twisted virtual handlebody-links and define
an integer-valued invariant $\Phi_{X}^{\mathbb{Z}}$ of twisted virtual 
handlebody-links. We provide example computations of the new invariants
and use them to distinguish some twisted virtual handlebody-links.
\end{abstract}

\parbox{5.25in} {\textsc{Keywords:} Twisted virtual handlebody-links, twisted virtual bikeigebras, counting invariants

\smallskip

}

\section{\large\textbf{Introduction}}\label{I}


 A \textit{spatial graph} is a graph embedded in $\mathbb{R}^3$, which we
may conceptualize as a ``knotted graph.'' A regular neighborhood of 
a trivalent spatial graph is a handlebody embedded in $\mathbb{R}^3$,
known as a \textit{handlebody-knot}. In \cite{A} Reidemeister moves for 
handlebody-knots are described and in \cite{AKMM} genus-two handlebody-knots 
with up to six crossings are classified.

In \cite{K}, virtual knot theory was introduced, with a geometric
interpretation of virtual knots as knotted curves in thickened orientable
surfaces up to stabilization moves described in \cite{CKS}. In \cite{B}, this 
geometric interpretation of virtual knots and link was extended to include
knotted curves in thickened non-orientable surfaces known as \textit{twisted 
virtual links}.

\textit{Bikei}, also known as \textit{involutory biquandles}, were considered
in \cite{AN}; twisted biracks and their special case, twisted bikei, were
considered in \cite{CN}. These algebraic structures define invariants of
unoriented twisted virtual links by counting colorings of the semiarcs
in a twisted virtual link diagram satisfying certain conditions at the
crossings and twist bars.

In this paper we consider \textit{twisted virtual handlebody-knots}, which are
handlebodies embedded in 3-manifolds with boundary of the form 
$\Sigma\times [0,1]$ for a possibly non-orientable surface $\Sigma$. We 
identify a list of Reidemeister moves for diagrams representing these objects 
and use these moves to define an algebraic structure we call 
\textit{twisted virtual bikeigebras}, generalizing algebraic structures such 
as those in \cite{CN, KK2, L}.
We then show that twisted virtual bikeigebra colorings of unoriented
twisted virtual handlebody-link diagrams are preserved by our moves in 
the sense that for every coloring before a move, there is a unique 
corresponding coloring after the move; hence, the set of such colorings
is an invariant of twisted handlebody-knots. 

The paper is organized as follows. In Section \ref{T} we introduce twisted 
virtual handlebody-knots. In Section \ref{B} we introduce twisted virtual
bikeigebras and use them to define computable invariants of twisted virtual 
handlebody-links with explicit examples.
We end in Section \ref{Q} with some questions for future work.

\section{\large\textbf{Twisted Virtual Handlebody-Knots}}\label{T}

We begin with a few basic definitions.

A \textit{handlebody-knot} or (\textit{handlebody-link} if there is more than
one component) is a handlebody of any genus embedded in $S^3$. In \cite{A}, 
handlebody-knots are shown to be equivalent via taking regular neighborhoods 
to equivalence classes of trivalent spatial graphs modulo a set of moves 
including the classical Reidemeister moves together with additional moves 
involving the interaction of trivalent vertices with classical crossings.

A \textit{twisted virtual link} is an ambient isotopy class of simple closed 
curves in an ambient space of the form $\Sigma\times [0, 1]$, where $\Sigma$ is 
a compact surface, up to stabilization of $\Sigma$. If $\Sigma$ is orientable, 
we have a \textit{virtual link} and if $\Sigma$ is non-orientable, we have a 
\textit{twisted virtual link.} In \cite{B,CKS,KK,K}, virtual links and  
twisted virtual links were show to be equivalent to equivalence classes of 
diagrams including, in addition to the usual classical crossings, 
\textit{virtual crossings} representing genus in $\Sigma$ and 
\textit{twist bars} representing cross-caps in $\Sigma$. 
Starting with a twisted virtual link diagram $D$, we build an abstract surface 
$\Sigma$ by ``skimming'' from $D$ as shown:
\[\includegraphics{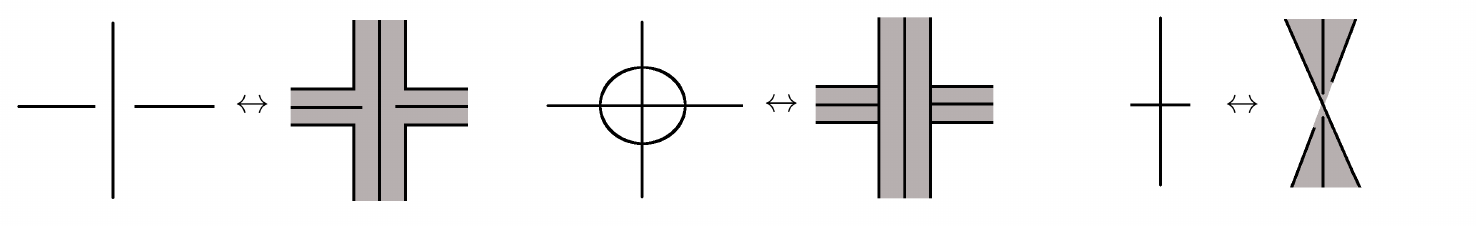}\]
We may think of our knot as the center line of a road with classical 
crossings representing intersections, virtual crossings representing bridges 
and twist bars representing the highway twisting upside down like a roller  
coaster. We cap off the resulting boundary circles with discs, possibly 
stabilizing with handles, to complete the abstract supporting surface $\Sigma$.

\begin{example}
The twisted virtual knot diagram below represents a knot in a thickened 
$T^2\#\mathbb{R}P^2$.
\[\includegraphics{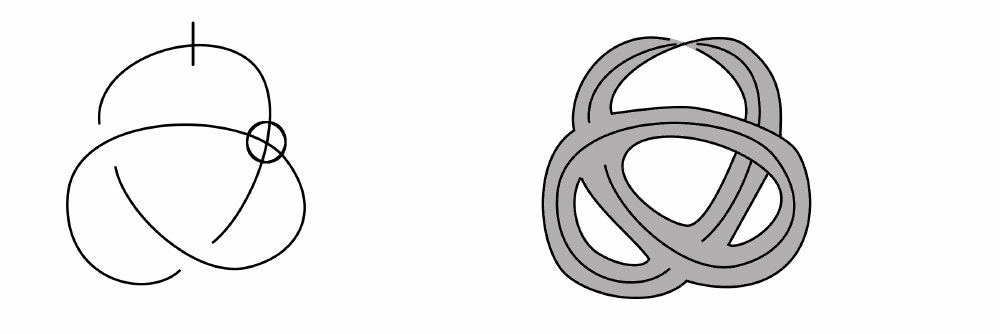}\]
\end{example}

We are interested in the case of twisted virtual handlebody-knots and 
twisted virtual handlebody-links. Geometrically, these are embeddings of 
handlebodies into ambient spaces of the form $\Sigma \times [0, 1]$
regarded up to ambient isotopy of $\Sigma\times [0,1]$ and stabilization 
of $\Sigma$. Considering the interaction of virtual crossings, twist bars and 
trivalent vertices, the reasoning in \cite{A} and \cite{B} yields two new 
additional moves as described below. 

\begin{remark}
One can consider the special case of virtual handlebody-links without allowing
twist bars; in this case, the analog of the Reidemeister IV move
replacing the classical crossing with a virtual crossing is treated as 
a forbidden move.
\end{remark}

\begin{definition}
\textit{Twisted virtual handlebody-knots} are equivalence classes of unoriented
trivalent spatial graph diagrams under the following rather lengthy
set of moves. First, equivalence classes of diagrams containing only classical 
crossings under the moves
\[\includegraphics{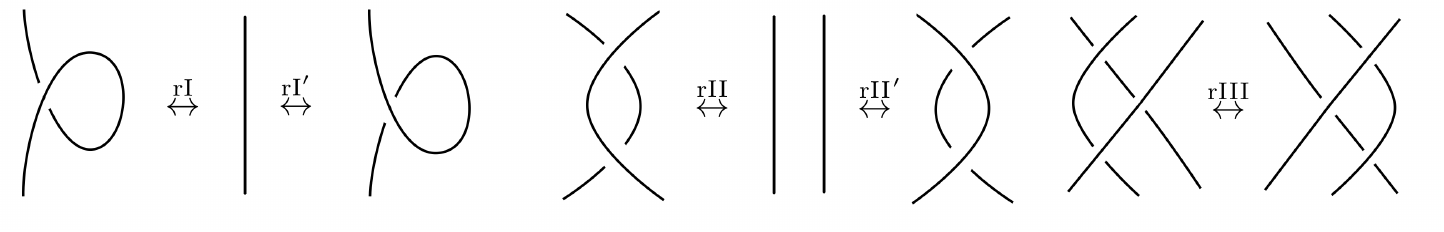}\]
form \textit{classical knots and links}. Including virtual crossings in our 
diagrams and allowing the moves
\[\includegraphics{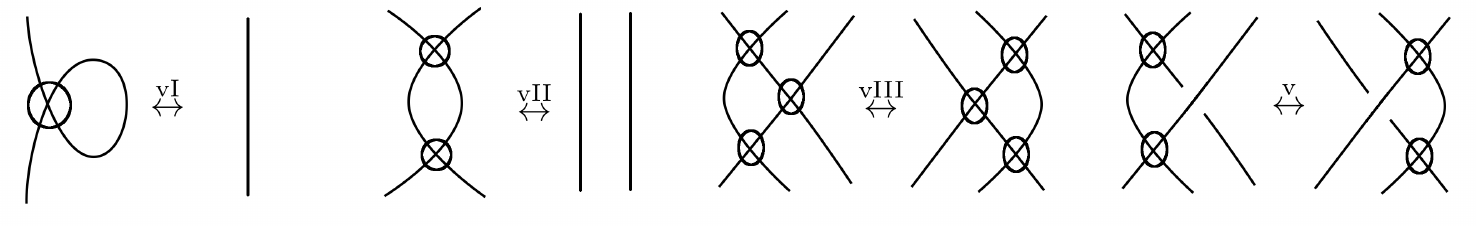}\]
yields \textit{virtual knots and links} \cite{K}. Allowing classical crossings 
and trivalent vertices together with moves rI,rII rIII and
\[\includegraphics{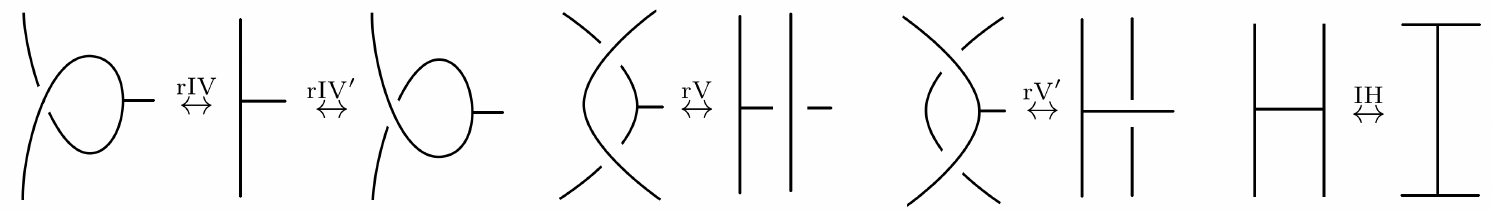}\]
yields \textit{handlebody-knots} \cite{A,AKMM}. Allowing classical crossings,
virtual crossings and twist bars with moves rI, rII, rIII, vI, vII, vIII and v 
and
\[\includegraphics{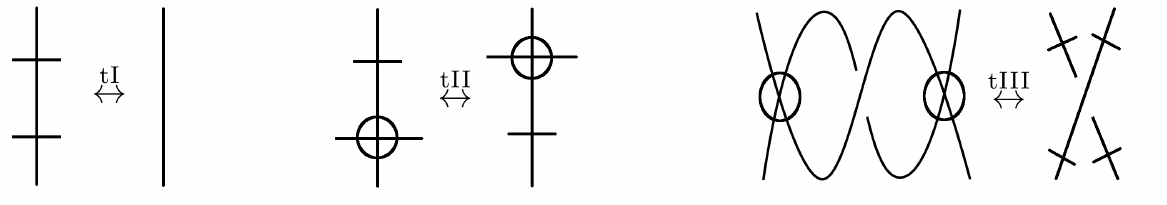}\]
yields \textit{twisted virtual knots and links} \cite{B}.
 
Since the interaction of trivalent vertices with classical crossings is 
considered in \cite{A} and the interaction of twists bars and virtual crossings
with classical crossings is considered in \cite{B,CKS,KK,K}, it remains
only to consider the interaction of twist bars and virtual crossings with
trivalent vertices. Thinking in terms of the skimming process, we have
the following two pictures
\[\includegraphics{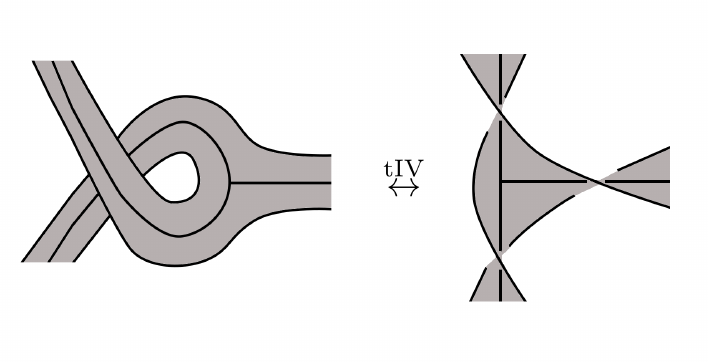}\includegraphics{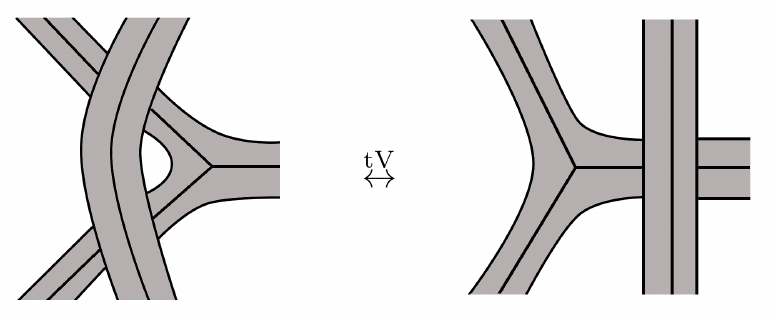}\]
which yield the moves tIV and tV below. Hence, allowing classical crossings, 
virtual crossings, trivalent vertices and twist bars with all previously 
listed moves, we need only two more moves, namely 
\[\includegraphics{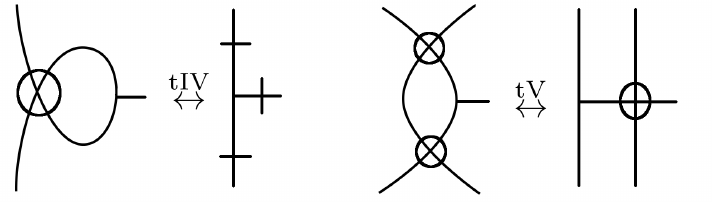}\]
to define \textit{twisted virtual handlebody-knots}.
\end{definition}

\section{\large\textbf{Twisted Virtual Bikeigebras}}\label{B}

\begin{definition}
Let $X$ be a set. A \textit{twisted virtual bikeigebra} structure on $X$
consists of three binary operations $x,y\mapsto x\utr y,x\otr y,x\vtr y$
and one partially defined operation $(x,y)\mapsto xy$
and an involution $T:X\to X$ satisfying the following axioms:
\begin{itemize}
\item[(i)] For all $x\in X$, we have $x\utr x=x\otr x\quad (ri.i)$,
\item[(ii)] For all $x,y\in X$ we have
\[\begin{array}{rcll}
(x\utr y)\utr y & = & x & (rii.i) \\
(x\otr y)\otr y & = & x & (rii.ii )\\
\end{array}\quad
\begin{array}{rcll}
x\utr(y\otr x) & = & x\utr y & (rii.iii)\\
x\otr(y\utr x) & = & x\otr y & (rii.iv)\\
\end{array}\quad
\begin{array}{rcll}
(x\vtr y)\vtr y & = & x & (vii.i)\\
x\vtr(y\vtr x) & = & x\vtr y & (vii.ii) \\
\end{array}
\]
\item[(iii)] For all $x,y,z\in X$, we have the classical and virtual
\textit{exchange laws}:
\[\begin{array}{rcll}
(x\utr y)\utr(z\utr y) & = & (x\utr z)\utr (y\otr z) & (riii.i) \\
(x\utr y)\otr(z\utr y) & = & (x\otr z)\utr (y\otr z) & (riii.ii) \\
(x\otr y)\otr(z\otr y) & = & (x\otr z)\otr (y\utr z) & (riii.iii) \\
\end{array}\quad
\begin{array}{rcll}
(x\utr y)\vtr(z\vtr y) & = & (x\vtr z)\utr (y\vtr z) & (viii.i)\\
(x\otr y)\vtr(z\vtr y) & = & (x\vtr z)\otr (y\vtr z) & (viii.ii)\\
(x\vtr y)\vtr(z\otr y) & = & (x\vtr z)\vtr (y\utr z) & (viii.iiii)\\
(x\vtr y)\vtr(z\vtr y) & = & (x\vtr z)\vtr (y\vtr z) & (v.i) \\
\end{array}
\]
\item[(iv)] For elements where the operation $(x,y)\mapsto xy$ is defined, 
the operation is associative and we have 
\[ 
 xy=z \iff  yz=x  \iff  zx=y,
\quad
\begin{array}{rcll}
xy & = & (y\otr x)(x\utr y) & (riv.i)\\
xy & = & (y\utr x)(x\otr y) & (riv.ii)\\
\end{array}
\]
and
\[
\begin{array}{rcll}
(xy)\utr z & = & (x\utr z)(y\utr (z\otr x))  & (rv.i) \\
(xy)\otr z & = & (x\otr z)(y\otr (z\utr x))  & (rv.ii) \\
(xy)\vtr z & = & (x\vtr z)(y\vtr (z\vtr x))  & (tv.i) \\
\end{array}\begin{array}{rcll}
x\utr(yz) & = & (x\utr y)\utr z & (rv.iii) \\
x\otr(yz) & = & (x\otr y)\otr z & (rv.iv) \\
x\vtr(yz) & = & (x\vtr y)\vtr z & (tv.ii) \\
\end{array}
\]
\item[(v)]
For all $x,y\in X$,
\[
\begin{array}{rcll}
T(x\vtr y) & = & T(x)\vtr y & (tii.i) \\
x\vtr y & = & x\vtr T(y) & (tii.ii)  \\
T(T(x)T(y)) & = & (y\vtr x)(x\vtr y) & (tiv.i) \\
T(T(x)\utr T(y)) & = & [(x\vtr y)\otr(y\vtr x)]\vtr[(y\vtr x)\utr(x\vtr y)] & (tiii.i) \\
T(T(x)\otr T(y)) & = & [(x\vtr y)\utr(y\vtr x)]\vtr[(y\vtr x)\otr(x\vtr y)] &(tiii.ii)\\
\end{array}
\]
\end{itemize}
\end{definition}

\begin{example}
Let $X$ be a set and define trivial bikei and virtual bikei operations on $X$,
so we have for all $x,y\in X$
\[x\utr y=x\otr y= x\vtr y =x.\]
Then the conditions of axioms (i), (ii) and  (iii) are satisfied automatically,
as are conditions (rv.i-iv) and (tv.i-ii); the rest of axiom (iv) requires that
the partially defined multiplication on $X$ be associative and commutative when 
defined, with \[xy=x\iff yz=x\iff zx=y.\] 
Writing $xy$ additively, this condition
says \[x+y=z\iff y+z=x\iff z+x=y;\]
we can satisfy this by making $X$ a $\mathbb{Z}_2$-module so that all three
equations become simply
\[x+y+z=0.\]
For example, the Klein 4-group $X=\mathbb{Z}_2\oplus\mathbb{Z}_2$ is a
virtual bikeigebra with $x\otr y=x\utr y=x\vtr y=x$ and $xy=x+y$.
Then any involution $T:X\to X$ satisfies (tii.i-ii) and (tiii.i-ii) 
automatically, and the final condition (tiv.i) requires that
\[T(T(x)+T(y))=y+x\]
or equivalently
\[T(x)+T(y)=T(x+y),\]
so the twist maps admitted by $X$ are the automorphisms of $X$.
\end{example}

\begin{example}\label{alex}
An \textit{Alexander bikei} is a module $X$ over 
$\mathbb{Z}[t^{\pm 1}, s^{\pm 1}]/(s^2-1, t^2-1, (1-s)(s-t))$
with bikei operations
\[\begin{array}{rcl}
x\utr y & = & tx+(s-t) y \\
x\otr y & = & sx
\end{array}.\]
Then if $X$ is a twisted virtual bikeigebra, axiom (rv.iv) requires
\[sx=x\otr(yz)=(x\otr y)\otr z = s^2x=x\]
so we must have $s=1$; 
then axiom (rv.iii) requires 
\[tx+(1-t)yz=x\utr(yz)=(x\utr y)\utr z
=t^2x+t(1-t)y+(1-t)z 
=x+t(1-t)y+(1-t)z 
\]
so we must also have $t=1$; hence any twisted virtual bikeigebra with
Alexander bikei operations must have trivial bikei operations
\[\begin{array}{rcl}
x\utr y & = & 1x+(1-1)y=x \\
x\otr y & = & 1x =x.
\end{array}\]
\end{example}

More generally, we can specify a twisted virtual bikeigebra with a 
block matrix encoding the operation tables of the three operations
$\utr,\otr,\vtr$, one partially defined operation $(x,y)\mapsto xy$, 
and one involution $T$.

\begin{example}
The set $X=\{1,2,3,4\}$ has twisted virtual bikeigebra structures including
\[
\begin{array}{r|rrrr}
\utr & 1 & 2 & 3 & 4 \\ \hline
1 & 2 & 2 & 1 & 1 \\
2 & 1 & 1 & 2 & 2 \\
3 & 3 & 3 & 3 & 3 \\
4 & 4 & 4 & 4 & 4 \\
\end{array}\quad
\begin{array}{r|rrrr}
\otr & 1 & 2 & 3 & 4 \\ \hline
1 & 2 & 2 & 1 & 1 \\
2 & 1 & 1 & 2 & 2 \\
3 & 3 & 3 & 3 & 3 \\
4 & 4 & 4 & 4 & 4 \\
\end{array}\quad
\begin{array}{r|rrrr}
\vtr & 1 & 2 & 3 & 4 \\ \hline
1 & 2 & 2 & 1 & 1 \\
2 & 1 & 1 & 2 & 2 \\
3 & 3 & 3 & 3 & 3 \\
4 & 4 & 4 & 4 & 4 \\
\end{array}\quad
\begin{array}{r|rrrr}
\cdot & 1 & 2 & 3 & 4 \\ \hline
1 & 3 & 4 & 1 & 2 \\
2 & 4 & 3 & 2 & 1 \\
3 & 1 & 2 & 3 & 4 \\
4 & 2 & 1 & 4 & 3 \\
\end{array}\quad
\begin{array}{cc}
x & T(x) \\ \hline
1 & 2 \\ 
2 & 1 \\
3 & 3 \\
4 & 4
\end{array}
\]
which we can write more compactly as a block matrix
\[\left[
\begin{array}{rrrr|rrrr|rrrr|rrrr|r}
2 & 2 & 1 & 1 & 2 & 2 & 1 & 1 & 2 & 2 & 1 & 1 & 3 & 4 & 1 & 2 & 2 \\
1 & 1 & 2 & 2 & 1 & 1 & 2 & 2 & 1 & 1 & 2 & 2 & 4 & 3 & 2 & 1 & 1 \\
3 & 3 & 3 & 3 & 3 & 3 & 3 & 3 & 3 & 3 & 3 & 3 & 1 & 2 & 3 & 4 & 3 \\
4 & 4 & 4 & 4 & 4 & 4 & 4 & 4 & 4 & 4 & 4 & 4 & 2 & 1 & 4 & 3 & 4 \\
\end{array}\right].\]
The same bikei can have several different virtual structures,
partial products and twist maps; for example, we also have twisted virtual 
bikeigebra
\[\left[
\begin{array}{rrrr|rrrr|rrrr|rrrr|r}
2 & 2 & 1 & 1 & 2 & 2 & 1 & 1 & 1 & 1 & 1 & 1 & - & - & - & - & 2 \\
1 & 1 & 2 & 2 & 1 & 1 & 2 & 2 & 2 & 2 & 2 & 2 & - & - & - & - & 1 \\
3 & 3 & 3 & 3 & 3 & 3 & 3 & 3 & 3 & 3 & 3 & 3 & - & - & 3 & 4 & 3 \\
4 & 4 & 4 & 4 & 4 & 4 & 4 & 4 & 4 & 4 & 4 & 4 & - & - & 4 & 3 & 4 \\
\end{array}\right]\]
where $-$ indicates an undefined product.
\end{example}

The twisted virtual bikeigebra axioms are motivated by the Reidemeister 
moves for unoriented twisted virtual handlebody-knots with the coloring rules
\[\includegraphics{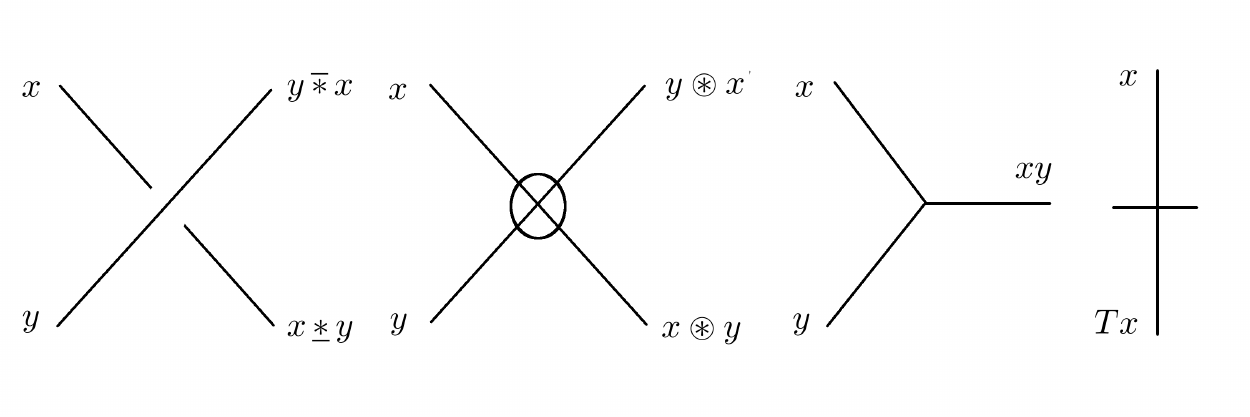}.\]

For instance, the move tV gives us axiom (tv.i):
\[\includegraphics{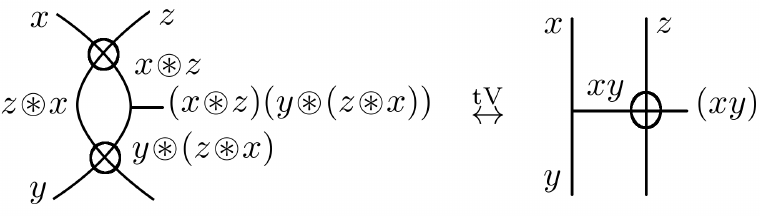}.\]

Analogously to other algebraic structures such as groups, quandles, biquandles
etc. \cite{EN}, we can define computable integer-valued invariants of twisted 
virtual handlebody-links by counting colorings of a diagram $D$ of a twisted 
virtual handlebody-link by elements of a twisted virtual bikeigebra $X$ 
satisfying the coloring rule above. By construction, we have the following 
theorem:

\begin{theorem}\label{thm:main}
The sets of colorings of twisted virtual handlebody-link diagram by a twisted
virtual bikeigebra before and after each Reidemeister move are in one-to-one 
correspondence.
\end{theorem}

This motivates the following definition:
\begin{definition}
Let $D$ be a twisted virtual handlebody-link diagram and $X$ a twisted virtual
bikeigebra. Then the cardinality of the set $\mathcal{C}(D,X)$ of colorings
of $D$ by $X$ is denoted $\Phi_X^{\mathbb{Z}}(D)=|\mathcal{C}(D,X)|$.
\end{definition}

Theorem \ref{thm:main} implies the following corollary:

\begin{corollary}
If $D$ and $D'$ are twisted virtual handlebody-links related by twisted 
virtual handlebody Reidemeister moves and $X$ is a twisted virtual bikeigebra, 
then $\Phi_X^{\mathbb{Z}}(D)=\Phi_X^{\mathbb{Z}}(D')$. In particular, 
$\Phi_X^{\mathbb{Z}}$ is an invariant of twisted virtual handlebody-links.
\end{corollary}

\begin{definition}
Let $K$ be a twisted virtual handlebody-link and $X$ a twisted virtual 
bikeigebra. If $\Phi_X^{\mathbb{Z}}(K)\ne 0$, we say $K$ is 
\textit{$X$-colorable}.
\end{definition}

\begin{remark}
We note that $X$-colorability is a two-valued invariant 
(``Yeah, bikeigebra-colorable by $X$'' or 
``No, not bikeigebra-colorable by $X$'') analogous to
the usual $p$-colorability of knots, with the distinction that 
a twisted virtual handlebody-link $L$ is  $X$-colorable is there is 
\textit{any} $X$-coloring of $L$; we do not need a notion of trivial 
vs. non-trivial colorings.
\end{remark}

\begin{example}\label{ex:1}
Let us compute the sets of colorings of the twisted virtual 
handlebody-knot $L_1$ below by the 
twisted virtual bikeigebra $X$ on the set $\{1,2\}$ given by the matrix 
\[\left[\begin{array}{rr|rr|rr|rr|r}
1& 1&  1& 1& 1& 1& 1 & - & 2 \\
2& 2&  2& 2& 2& 2& - & 2 & 1 \\
\end{array}\right].\]
\[\includegraphics{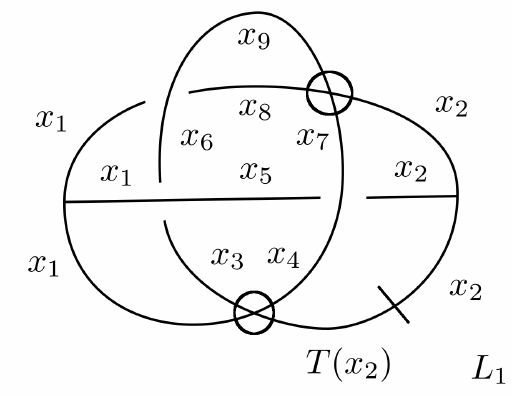}\]
The system of coloring equations is
\[\begin{array}{rcl}
x_3 & = & T(x_1)\vtr x_1 \\
x_4 & = & x_1 \vtr T(x_2) \\
x_5 & = & x_1 \otr x_3 \\
x_5 & = & x_2 \utr x_4 \\
x_6 & = & x_3 \utr x_1 \\
x_7 & = & x_4 \otr x_2 \\ 
x_8 & = & x_1 \utr x_6 \\
x_8 & = & x_2 \vtr x_7 \\
x_9 & = & x_6 \otr x_1 \\
x_9 & = & x_7 \vtr x_2  
\end{array} \leftrightarrow
\begin{array}{rcll}
x_1 \otr (T(x_2)\vtr x_1)  & = & x_2 \utr (x_1 \vtr T(x_2)) & (i) \\
x_1 \utr ((T(x_2)\vtr x_1) \utr x_1)  & = & x_2 \vtr ((x_1 \vtr T(x_2)) \otr x_2)  & (ii)\\
((T(x_2)\vtr x_1) \utr x_1) \otr x_1 
 & = & ((x_1 \vtr T(x_2)) \otr x_2) \vtr x_2 & (iii)
\end{array} 
\]
We can then verify for each of the four possible colorings that the system is
not satisfied:
\[\begin{array}{rr|rrr}
x_1 & x_2 & (i) & (ii) & (iii) \\ \hline
1 & 1 & 1=1 & 1\ne 1 & 2\ne 1 \\
1 & 2 & 1\ne 2 & 1\ne 2 & 2=2 \\
2 & 1 & 2\ne 1 & 2\ne 1 & 1=1  \\
2 & 2 & 2\ne 2 & 2=2 & 1\ne 2
\end{array}
\]
so the twisted virtual handlebody-knot $L_1$ is not $X$-colorable. 

On the other hand, the twisted virtual handlebody-knot below \textit{is}
$X$-colorable in two ways:
\[\includegraphics{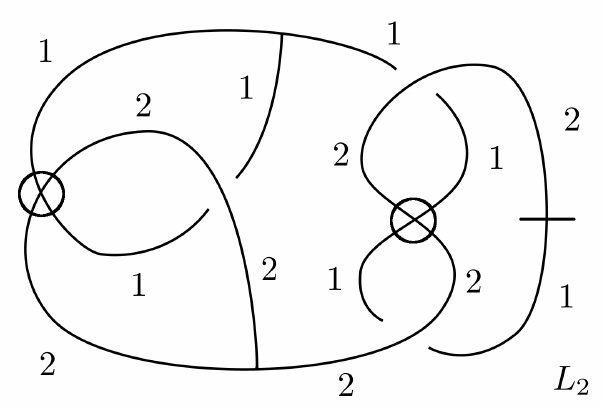}\includegraphics{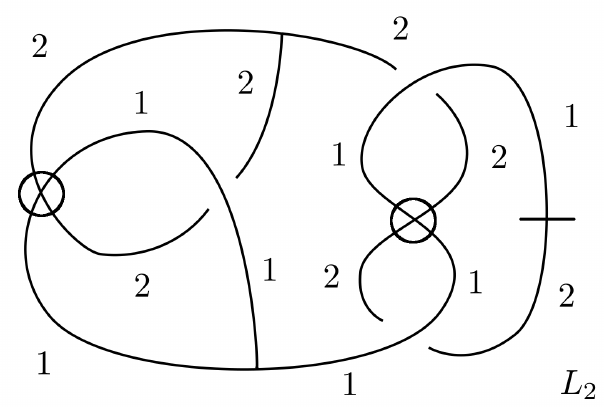}\]
and hence $L_1$ is distinguished from $L_2$ by $X$-colorability.
\end{example}

\begin{example}
Let $X$ be the twisted virtual bikeigebra on the set $\{1,2,3,4\}$ given
by the operation matrix
\[\left[\begin{array}{rrrr|rrrr|rrrr|rrrr|r}
1& 1& 1& 2& 1& 1& 2& 1& 1& 1& 2& 1 & 1 & - & - & - & 2 \\
2& 2& 2& 1& 2& 2& 1& 2& 2& 2& 1& 2 & - & 2 & - & - & 1 \\
3& 3& 4& 4& 3& 3& 4& 4& 3& 3& 3& 3 & - & - & - & - & 4 \\
4& 4& 3& 3& 4& 4& 3& 3& 4& 4& 4& 4 & - & - & - & - & 3 \\
\end{array}\right].\] Then the 
twisted virtual handlebody-link $L_3$ below is not $X$-colorable,
\[\includegraphics{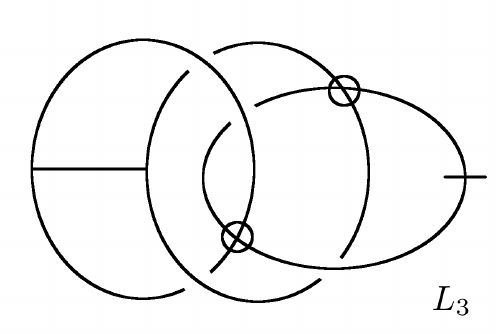}\]
but switching the position of the twist bar yields a
twisted virtual handlebody-link $L_4$ with two $X$-colorings.
\[\includegraphics{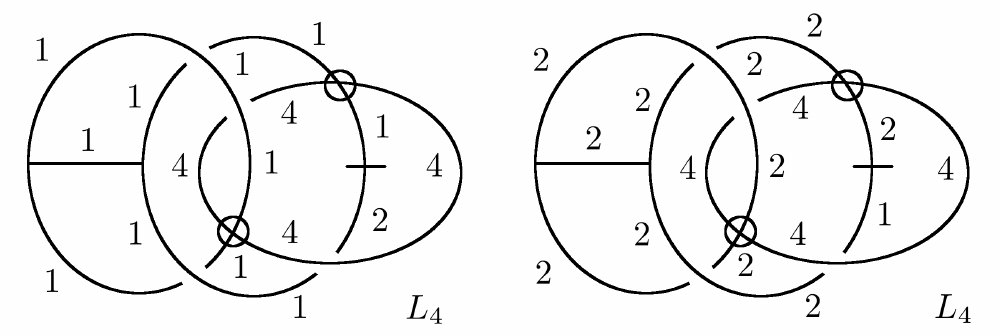}\]
\end{example}

When two twisted virtual handlebody-links are both colorable by the same
twisted virtual bikeigebra, we can use the number of colorings to distinguish 
the links.

\begin{example}
Consider the two twisted virtual handlebody-knots $L$ and $L'$ below, both 
of which are handlebody-knots of genus 2 represented by diagrams 
with four classical crossings, two 
virtual crossings and one twist bar. They are distinguished by their counting 
invariants with respect to the twisted virtual bikeigebra $X$ given by the 
operation matrix
\[\scalebox{0.9}{$
\left[\begin{array}{rrrrrrrr|rrrrrrrr|rrrrrrrr|rrrrrrrr|r}
1&1&1&1&1&1&1&1&1&1&1&1& 1& 1& 1& 1&1& 1& 1& 1& 1& 1& 1& 1& 1&2&3&4&5&6&7&8&1\\
2&6&6&2&2&6&6&2&2&6&2&6& 2& 6& 2& 6&2& 2& 2& 2& 2& 2& 2& 2& 2&1&4&3&6&5&8&7&2\\
3&3&7&7&3&3&7&7&3&7&7& 3& 3& 7& 7& 3&3& 3& 7& 7& 3& 3& 7& 7& 3&4&1&2&7&8&5&6&4\\
4&8&4&8&4&8&4&8&4&4&8& 8& 4& 4& 8& 8&4& 4& 8& 8& 4& 4& 8& 8& 4&3&2&1&8&7&6&5&3\\
5&5&5&5&5&5&5&5&5&5&5& 5& 5& 5& 5& 5&5& 5& 5& 5& 5& 5& 5& 5& 5&6&7&8&1&2&3&4&5\\
6&2&2&6&6&2&2&6&6&2&6& 2& 6& 2& 6& 2&6& 6& 6& 6& 6& 6& 6& 6& 6&5&8&7&2&1&4&3&6\\
7&7&3&3&7&7&3&3&7&3&3& 7& 7& 3& 3& 7&7& 7& 3& 3& 7& 7& 3& 3& 7&8&5&6&3&4&1&2&8\\
8&4&8&4&8&4&8&4&8&8&4& 4& 8& 8& 4& 4&8& 8& 4& 4& 8& 8& 4& 4&8&7&6&5&4&3&2&1&7\\
\end{array}\right]$}\]
with 40 and 64 colorings respectively as computed by our \texttt{Python} code.
\[\includegraphics{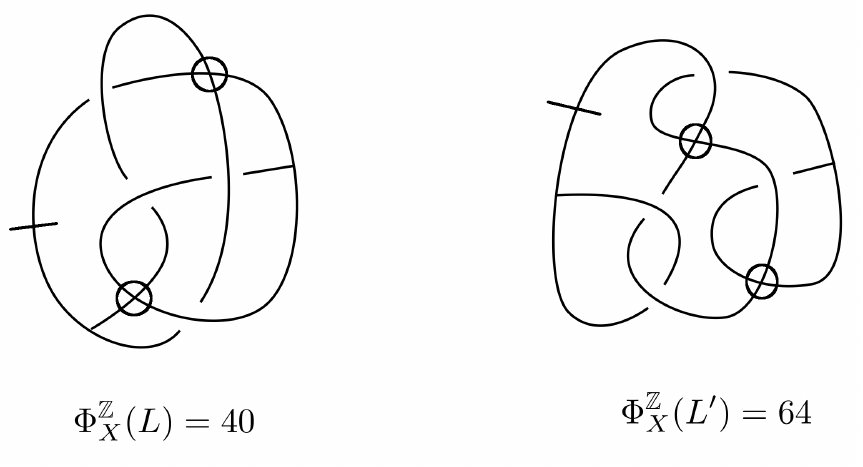}\]
\end{example}

\section{\large\textbf{Questions}}\label{Q}

We end with some questions for future research. 

We have only looked at the unoriented case of twisted virtual handlebody-links 
and their associated algebraic structure; the oriented case should define a
generalized algebraic structure we could call \textit{twisted virtual 
biqualgebras}. What are the appropriate axioms and examples of such structures?

We are very interested in examples of families of twisted virtual bikeigebra
structures defined in terms of groups, modules, matrix algebras, and any
other algebraic structures whose additional properties can be used to simplify
coloring calculations.

As with all counting invariants, we ask what enhancements of the twisted
virtual bikeigebra counting invariant are possible. Cocycle invariants, 
structure enhancements, twisted virtual bikeigebra module enhancements and 
many more possibilities remain to be investigated.

\bibliography{sn-yz}{}
\bibliographystyle{abbrv}

\bigskip

\noindent
\textsc{Department of Mathematical Sciences \\
Claremont McKenna College \\
850 Columbia Ave. \\
Claremont, CA 91711}

\end{document}